\documentclass[12pt]{amsart}
 \usepackage{graphicx}
\usepackage{amsmath,amsthm}
\usepackage[small,nohug, heads=LaTeX]{diagrams}
\diagramstyle[labelstyle=\scriptstyle]

\theoremstyle{plain}
\newtheorem{thm}{Theorem}[section]
\newtheorem{lemma}[thm]{Lemma}

\theoremstyle{definition}
\newtheorem{defn}[thm]{Definition}
\newtheorem{rmk}[thm]{Remark}

\newcommand{\pE}{\hbox{$^\pi\kern-2pt E$}}
\newcommand{\hQ}{\hbox{$\hat Q$}}
\newcommand{\phQ}{\hbox{$ '{\hat Q}$}}

\newtheorem{ex}[thm]{Example}

\newcommand{\K}{{\mathcal K}}

\newcommand{\la}{\langle}
\newcommand{\ra}{\rangle}

\newcommand{\onto}{{\protect \rightarrow\!\!\!\!\!\rightarrow}}

\newcounter{letter}
\renewcommand{\theletter}{\rom{(}\alph{letter}\rom{)}}

\newcounter{rnum}
\renewcommand{\thernum}{\rom{(}\roman{rnum}\rom{)}}

\begin{document}


\title{PBW-Deformation Theory and Regular Central Extensions}


\author[Cassidy and Shelton]{ }

\maketitle

\begin{center}

\vskip-.2in Thomas Cassidy$^\dagger$ and Brad Shelton$^\ddagger$\\
\bigskip

$^\dagger$Department of Mathematics\\ Bucknell University\\
Lewisburg, Pennsylvania  17837
\\ \ \\

   $^\ddagger$Department of Mathematics\\ University of Oregon\\
Eugene, Oregon  97403-1222
\end{center}

\begin{abstract}
\baselineskip15pt
A deformation $U$, of a graded $K$-algebra $A$ is said to be of PBW type if $grU$ is 
$A$.  It has been shown for Koszul and $N$-Koszul algebras that the deformation is PBW if and only if the relations of $U$ satisfy a Jacobi type condition.  In particular, for these algebras the determination of the PBW property is a {\it finite} and explicitly determined linear algebra problem.  We extend these results to an arbitrary graded $K$-algebra, using the notion of central extensions of algebras and a homological constant attached to $A$ which we call the {\it complexity} of $A$.  
\end{abstract}

\setcounter{page}{1}

\thispagestyle{empty}

\baselineskip18pt

\section{Introduction}

We are interested in certain deformations of graded algebras $A$ over a field $K$ obtained by altering the defining relations of $A$.    For reasons  explained below, these deformations are called PBW-deformations.  An example is  the universal enveloping algebra $U \frak g$ of a Lie algebra $\frak g$, which is a deformation of the symmetric algebra $\mathbb S\frak g$.   In general, determining when a given deformation $U$ of $A$ is a PBW-deformation  is a very difficult problem.  However it has been shown by  Polishchuk and Positselsky \cite{PP} and Braverman and Gaitsgory \cite{BG}  that for a Koszul algebra $A$ it suffices to verify a comparatively simple condition on the defining relations of $U$.  In the case $A=U\frak g$ this condition coincides with the Jacobi identity for the Lie bracket on $\frak g$, and so we call it the Jacobi condition.

   Motivated by the cubic Artin-Schelter (AS) regular algebras,  Berger \cite{B} introduced the notion of $N$-Koszul algebras as a generalization of Koszul algebras. A graded algebra $A$ is $N$-Koszul if all of its relations are in degree $N$ and its Yoneda Ext-algebra is generated by $Ext^1(K,K)$ and $Ext^2(K,K)$.  Two recent papers (\cite{BGZ} and \cite{FV})  
  have shown that a Jacobi type condition holds for $N$-Koszul algebras.   In particular, Fl\o ystad and Vatne \cite{FV} explicitly calculate the PBW-deformations of the cubic AS regular algebras.    Remarkably, the stipulations in their list of PBW-deformations are identical to those given previously in \cite{C1} for  central extensions of cubic AS regular algebras to be  AS regular.  Our work here was motivated by the desire to explain this phenomenon.

Throughout this paper the term {\it $K$-algebra} is used to denote an algebra over the field $K$ that is finitely generated and has a finitely generated set of relations, i.e.~an algebra of the form $$A=K\langle x_1\ldots x_n\rangle/\langle r_1,\ldots, r_m\rangle.$$ We may and will assume that none of the relations $r_i$ are linear. Such an algebra inherits a filtration, $F^k(A)$, from the standard filtration on the free algebra $K\langle x_1, \ldots ,x_n\rangle$.  

The free algebra $K\langle x_1\ldots ,x_n\rangle$ also carries a standard grading.  If the relations $r_1,\ldots,r_m$ are all homogeneous,  then $A$ inherits a grading from the free algebra.  We denote the component of degree $n$ by $A_n$ and we use the phrase {\it graded $K$-algebra} to describe this situation.  We  will assume that $r_1,\ldots,r_m$ is a minimal set of relations for $A$.  When $A$ is graded it is also augmented by $A \to A/A_+=K$.  

Let $A$ be a graded $K$-algebra as above, with relations $R=\{r_1,...,r_m\}$. Let 
$T$ be the free algebra $K\langle x_1,\ldots, x_n\rangle$.   By a {\it deformation} of $A$ we mean a $K$-algebra $U = T/\langle  r_1+l_1, \ldots ,r_m+l_m\rangle$, where 
$l_1,\ldots,l_m$ are (not necessarily homogenous) elements of $T$ such that $deg(l_i)<deg(r_i)$ for all $i$.  The algebra $U$ is not, in general, graded, but it has an associated graded algebra $gr(U) = \oplus gr(U)_k$, where 
$gr(U)_k = F^k(U)/F^{k-1}(U)$.  
 
 The classical Poincar\'e-Birkhoff-Witt Theorem   
 says that for  a Lie algebra $\frak g$ the associated graded ring $gr(U\frak g)$ is isomorphic to the symmetric algebra  $\mathbb S\frak g$.  This observation inspired the following definition (see \cite{PP} or \cite{BG}).

\begin{defn}  The non-graded deformation $U$, of the graded $K$-algebra $A$, is said to be a PBW-deformation if $grU$ is isomorphic to $A$.
\end{defn}
 PBW-deformations have been studied recently in \cite{EG}, \cite{GG} and \cite {CBH}. Deciding whether or not a given deformation is PBW appears to be an infinte linear algebra problem.  The primary goal of this paper is to determine the conditions under which this can be reduced to a {\it finite} linear calculation and to ascertain precisely what that linear problem entails.  We answer these questions by translating them into questions about a larger graded ring.  

Let $T[z] = T\otimes_K K[z]$ be the polynomial extension of $T$.  This is naturally graded by the K\"unneth formula.  Define a function $h:T \to T[z]$ as follows.  If $f = \sum_{i=0}^p f_i$, with $f_i \in T_i$ and $f_p \ne 0$, then $h(f) = \sum_{i=0}^p z^{p-i}f_i$.  We call $h(f)$ the {\it homogenization} of $f$ in $T[z]$.

Let $A$ be a graded $K$-algebra and $U$ a deformation of $A$, as above.  The {\it central extension} of the graded $K$-algebra $A$ {\it associated} to $U$ is the algebra 
$D:=T[z]/\langle h(r_1+l_1),\ldots , h(r_m+l_m)\rangle$.  This algebra is a graded $K$-algebra with one more generator than $A$ and  $n$ more relations (the commutation relations for $z$.)   

\begin{defn}  The central extension $D$ of $A$ associated to $U$ is said to be a {\it regular central extension} if the element $z$ in $D$ is not a zero divisor.  The central extension is {\it regular to degree $p$} if $Ann_D(z^n)_k = 0$ for all $k\le p$ and all $n\ge 1$.    
\end{defn}

We relate the notion of PBW-deformation directly to the notion of regular central extension through the following theorem.  

\begin{thm}\label{standard}  The deformation $U$ of the graded algebra $A$ is a PBW-deformation if and only if the associated central extension $D$ is a regular central extension.
\end{thm}

It is very simple to show that the Jacobi condition of \cite{BG} or \cite {PP} for deformations of Koszul algebras is equivalent to the statement that $z$ is not a zero divisor on $D_2$.  Likewise the Jacobi condition from \cite{BGZ} or \cite{FV} for  $N$-Koszul algebras  can be interpreted as the statement that $z$ is not a zero divisor on 
$D_{N}$ (and hence on $D_{\le N})$.  In both cases, the Jacobi condition assures that the deformation is PBW.  In light of the theorem above, this suggests that in general one may be able to determine the regularity of the  element $z$ from knowing its regularity only in some  specific low degree.  This is exactly the case and the appropriate degree is given by the following definition.

\begin{defn}\label{complexity}   The {\it complexity} of the graded $K$-algebra $A$ is 
$$c(A) = sup\{n\,|\, Ext^{3,n}(K,K) \ne 0\} - 1$$
if the global dimension of $A$ is at least 3.  For global dimension less than 3 we set $c(A)=0$.  

\end{defn}

This homological invariant of the algebra $A$ provides the precise measurement required to bound the size of the linear algebra problem attached to the PBW property.
We can now state our main theorem.  

\begin{thm}\label{main} Let $A$ be a graded $K$-algebra of finite complexity $c=c(A)$.  Let $D$ be a central extension of $A$ by the variable $z$.  Then $z$ is regular in $D$ if and only if $z$ is not a zero divisor on $D_{\le c}$.  
\end{thm}

Since the complexity of a Koszul algebra is 0 or 2,  and the complexity of an $N$-Koszul algebra is 0 or $N$, the PBW deformation results cited earlier in these two cases become special cases of \ref{main}.  On the other hand, a quadratic algebra can easily have complexity 2 without being Koszul, so even in the quadratic setting \ref{main} is seen to strengthen the classical result.  Moreover, the theorem has the advantage of dealing with algebras with relations in mixed degrees and does not in any way rely on homological purity.  

Theorem \ref{main} is a consequence of a more complete statement given in section 3 as Theorem \ref{main2}.  This Theorem also gives a simple criterion for regularity of $z$ through a single matrix equation involving the minimal projective resolution of the trivial module $K$,
$$\pi_D(M_3 f_2 + f_3 M_1) = 0.$$
This is probably the most compact and easily verified condition with which to check regularity.  We refer the reader to section 3 for the necessary notation.

It may be convenient to have the regularity criterion of Theorem \ref{main} reexpressed in terms of the original deformation $U$.  This is accomplished by translating the hypothesis: $z$ is not a zero divisor on $D_{\le c}$, back into a set of conditions on the relations of $U$.  We refer the reader to section 4 for the precise formulation of this translation, which we call a Jacobi condition.

Section 2 contains the proof of Theorem \ref{standard} in a slightly stronger form.  Section 3 introduces the necessary constructions and homological algebra required to prove Theorem \ref{main2}, and as a corollary Theorem \ref{main}.  In section 4 we give the Jacobi condition which solves the orginal problem and in section 5 we present a small number of interesting and motivating examples. 

Finally, we want to note the relationship of this paper to another line of inquiry.  The following lemma is easy to prove.

\begin{lemma}\label{K2}  If $Ext_A^3(K,K)$ is algebraically generated  by $Ext_A^1(K,K)$ and $Ext_A^2(K,K)$, then the complexity of $A$ is no greater than the maximum of the degrees of the relations of $A$.  
\end{lemma}

In other words, having $Ext_A(K,K)$ generated in low cohomological degree is an effective way to limit the size of $c(A)$. 
 Green et.~al.~\cite{GMMZ} have shown $N$-Koszul algebras satisfy the hypothesis of the lemma, and it is well known that  the class of Koszul algebras is precisely the class of graded algebras for which $Ext_A(K,K)$ is generated by $Ext_A^1(K,K)$.  The lemma, together with theorem \ref{main} leads us to make the following definition.

\begin{defn}  We say that the graded $K$-algebra $A$ is $\K_2$ if $Ext_A(K,K)$ is generated as an algebra by $Ext_A^1(K,K)$ and $Ext_A^2(K,K)$.  
\end{defn}

The class of $\K_2$ algebras will be studied in depth in the companion paper \cite{CS-K2}.
 
\section{Proof of Theorem \ref{standard}}
 
Theorem \ref{standard} will be an immediate consequence of the following stronger statement.  Throughout we have a fixed graded $K$-algebra $A$, a deformation $U$ and the associated central extension $D$. 

Let $R=\{r_1,\ldots,r_m\}$ be the set of relations of $A$ and 
$P=\{r_1+l_1,\ldots ,r_m+l_m\}$ the set of relations of $U$.  The set of relations of $D$ in $T[z]$ is then $h(P)$.  We note: $\langle h(P)\rangle \not= h(\langle P\rangle)$ in $T[z]$.
Let  $\pi:T\to grU$ be the natural epimorphism.  Since $\pi(R)=0$, we may factor $\pi$ through $A$ to obtain 
$\Phi:A \to gr(U)$ given by $\Phi(x+\langle R\rangle) = \pi(x)$.  This graded $K$-algebra homomorphism is always surjective. 

Let $\pi_D:T[z]\to D$ be the canonical epimorphism.  We consider $T$ to be a subalgebra of $T[z]$. 
Let $\phi_1, \phi_0:T[z]\to T$ be the maps which evaluate at $z=1$ and $z=0$ respectively.  These maps and those used in the next section are summarized in the following  diagram.

\begin{diagram}
&&T[z]& \pile{ \rOnto^{\ \ \phi_0\ \ \ }\\ \rOnto_{\ \ \phi_1\ \ \ }} &T & & &\\
&&\dOnto^{\pi_D}& & \dOnto(0,3)^{\pi_A} & \rdOnto(2,4)^\pi & \\
&&D& & & \\
  &\ldOnto^{} & & \rdOnto^{\phi}  && \\
 U&& &&A &  \rOnto_{\Phi\ \ }&grU
\end{diagram}

 \begin{thm} The map $\Phi:A \to grU$ is injective in degrees $1\le i\le p$  if and only the central extension $D$ is regular to degree $p$. 
\end{thm}

\begin{proof}  

Assume first $D$ is not regular to degree $p$.  This means that the set 
$$S = \{ x\in T[z]_{\le p}\setminus \langle h(P)\rangle\,|\, 
	z^kx \in \langle h(P)\rangle\, \hbox{\rm for some}\,k>0\}$$ is not empty.   
We choose $x\in S$ of minimal degree and choose the smallest $k$ for which $z^kx \in   \langle h(P)\rangle$.  We may then write $x=x'+zy$ for some $x'\in T$ and $y\in T[z]$ with $deg(y)<deg(x)$.  Since $z^k(x'+zy)\in \langle h(P)\rangle$, $\phi_1(z^k(x'+zy)) = x'+\phi_1(y)\in \langle P\rangle$.  It follows that 
$\pi(x')$ is 0 in $gr(U)$ and hence $x'+\langle R \rangle$ is in the kernel of $\Phi$.  

We claim that $x'$ is not in $\langle R\rangle$, which will prove that $\Phi$ is not injective on $A_{\le p}$.  Suppose to the contrary that 
$x'\in  \langle R\rangle$.  We may choose a representation of $x'$ as
$x' = \sum_{i=1}^m \sum_j u_{i,j}r_i w_{i,j}$.   Set 
$\tilde x = \sum_{i=1}^m \sum_j u_{i,j}h(r_i+l_i) w_{i,j} \in \langle h(P)\rangle$.  Since 
$h(r_i + l_i) - r_i$ is a multiple of $z$ (or zero) for every $i$, we see that
$\tilde x - x' = z y'$ for some $y'\in T[z]$.  We then have 
$z^k(x - \tilde x) = z^k(z(y-y')) \in \langle h(P)\rangle$.  By the minimality of $x$, we conclude $y-y'=0$.  But then $x = x' + zy' = \tilde x \in \langle h(P)\rangle$, which is a contradiction.     
  
For the converse, assume that $D$ is regular to degree $p$.   Choose a homogeneous $x\in T_{\le p}$ such that $x+\langle R \rangle$ is in the kernel of $\Phi$.  Then there is some (not necessarily homogeneous) $y\in T$ of lower degree than $x$, such that $x+y$ is in $\langle P\rangle$.   Let $k=deg(x)-deg(y)>0$.  Since 
$\phi_1^{-1}(\{x+y\})\cap T[z]_{hom}=\{z^j(x+z^kh(y))\,|\, j\ge 0\}$ there exists some non-negative integer $j$ such that $z^j(x+z^kh(y))$ is in 
$\langle h(P)\rangle$.  Since $z^j$ is not a zero divisor on $D_{\le p}$, this means 
$x+z^kh(y)$ is in $\langle h(P)\rangle$. Since $k>0$ we then get
$\phi_0(x+z^ky)=x\in \langle R \rangle$.  Hence $\Phi$ is injective, as required. 
\end{proof}

\section{Proof of Theorem \ref{main}}

Throughout this section, let $A$ and $D$ be as in the previous sections, that is $A$ is a graded $K$-algebra and $D$ a central extension of $A= D/Dz$, with $z\in D_1$ central in $D$.  Let $\phi:D\to A$ be the associated graded ring epimorphism (evaluation at 
$z=0$).  

We extend $\phi$ to the functor $\phi:$ Mod-$D \to $Mod-$A$ via $M \mapsto M/zM$.  We also use $\phi$ to denote the epimorphism $\phi:T[z]\to T$ given by evaluation at $z=0$. 

For any graded vector space $V=\oplus_k V_k$, let $V(j)$ be the same vector space with the shifted grading $V(j)_k=V_{j+k}$.  We will be dealing with graded-free resolutions of $A$ and $D$ modules, with generators in many different degrees, so it is convenient to introduce the notation
$A(j_1,\ldots,j_p) = \oplus_{i=1}^p A(j_i)$.

Let  $(Q, \delta)$ be a fixed minimal projective resolution of the trivial left $A$-module $A/A_+ = K$
$$\to Q^3 \buildrel \delta_3 \over \to Q^2 \buildrel \delta_2\over \to Q^1 \buildrel \delta_1\over\to Q^0 \to K\to 0.$$
   The modules $Q^k$ are graded-free modules, but may not be finitely generated.  However, the number of generators in any given degree must be finite since $A$ is locally finite as a graded vector space. Recall that $\pi_A:T\to A$ is the canonical epimorphism.  We choose a graded basis for each $Q^k$ and we choose a matrix 
$M_k$, with homogeneous entries from $T$, such that right multiplication by $\pi_A(M_k)$ represents $\delta_k:Q^k \to Q^{k-1}$.   Since the resolution $Q$ is minimal, we may assume that $Q^0=A$, 
$Q^1=A(-1,\ldots ,-1)$ and the matrix $M_1$ is the column $(x_1\,x_2\,\ldots\,x_n)^t$.  Similarly we may assume the finite matrix $M_2$ is chosen so that $M_2M_1 = (r_1,\ldots,r_m)^t$.  In general, the matrix 
$M_n$, for $n>2$, may not be finite, but local finiteness assures that it will have at most finitely many nonzero entries in every row and at most finitely many entries of a given degree in any column.  

Since the resolution $Q$ consists of graded-free modules, we may express each $Q^n$ as:
 $$Q^n = A(m_{1,n},m_{2,n},\ldots,m_{t_n,n}).$$
For $n>2$, $t_n$ may be infinite.  We will abuse notation slightly by treating $t_n$ as if it is finite, but the reader should have no trouble making the appropriate adjustment in the constructions that follow. Note that the shifts $m_{i,n}$ are all nonpositive, indeed $m_{i,n}\le -n$ for all $i$. The matrix $M_n$ is a $t_n$ by $t_{n-1}$-matrix, whose $i$-$j$ entry is a homogeneous element of $T$ of degree  $m_{i,n} - m_{j,n-1}$.  
 
We record the following obvious lemma.

\begin{lemma}\label{matrixcomplexity}  If the global dimension of $A$ is at least 3, then the complexity $c(A)$ is the supremum of the degrees of the entries of the matrix $M_3M_2$.  In particular, $c(A)$ is finite if and only if $M_3$ is a finite matrix.
\end{lemma}
 
We now turn to building sequences of graded-free modules for the central extension 
$D$.  Let $\hat Q^0 = D$ and for $n\ge 1$ put
 $$\hat Q^{n} = D(m_{1,n},m_{2,n},\ldots,m_{t_n,n}) \oplus D(m_{1,n-1},m_{2,n-1},\ldots,m_{t_{n-1},n-1})(-1).$$
We exhort the reader to remark on the extra degree shift in the second term of this definition. Let $I_t$ denote the $t$ by $t$ identity matrix.  For each $n> 1$, let $f_n$ be a $t_n$ by $t_{n-2}$ matrix with homogeneous entries from $T[z]$ which satisfies:
$$\pi_D(M_n M_{n-1} -  (-1)^{n-1}z f_n)=0.$$
Such $f_n$ exists because the entries of $M_nM_{n-1}$ are all in $\langle R\rangle$.
Next, define matrices $\hat M_n$, with homogeneous entries in $T[z]$, by
$$\hat M_1 = \left({\begin{array}{c}M_1\\ z \end{array}}\right)$$
and for $n\ge 2$
$$\hat M_n = \left({\begin{array}{cc} M_n & f_n \\ (-1)^{n-1}	zI_{t_{n-1}} & M_{n-1} \end{array}}\right).$$
These matrices may well be infinite, but they have the same finiteness conditions on rows and columns as the matrices $M_n$.  Finally, we define $\hat \delta_n:\hat Q^n \to \hat Q^{n-1}$ to be right multiplication by the matrix $\pi_D(\hat M_n)$ (and $\hat\delta_0=0$). 

Set $\hQ = \oplus_n \hat Q^n$ and $\hat\delta = \sum_n\hat \delta_n$. 
 
\begin{rmk}\label{corner}  There is no assurance that $(\hat Q, \hat \delta)$ is a complex. Indeed, determining when it is a complex is our main concern.  We do have, by construction, $\hat\delta_1 \hat\delta_2 = 0$ since $\pi_D(\hat M_2\hat M_1)=0$.  Hence  $\hQ^2 \to \hQ^1 \to \hQ^0\to 0$, at least, is a complex.  It is not difficult to see that the homology of this complex is $H^0(\hQ) = K$ and $H^1(\hQ) = 0$.

For $n>2$, $ \hat \delta_{n-1} \hat\delta_{n}$ is given by right multiplication by the matrix $\pi_D(\hat M_n \hat M_{n-1})$, and by the centrality of $z$, we see that this matrix is 
$$ 	\left(\begin{array}{cc}
	0 & \pi_D(M_nf_{n-1} + f_n M_{n-2}) \\
	0 & 0
		\end{array}\right).
		 $$
\end{rmk}

\begin{lemma}\label{complex1} Suppose that $(\hat Q, \hat \delta)$ is a complex.  Then the $A$-module complex
$(\phi(\hat Q),\phi(\hat \delta))$ has the homology:
$H^0(\phi(\hat Q)) = K$, $H^1(\phi(\hat Q)) = K(-1)$ and $H^p(\phi(\hat Q)) = 0$ for 
$p>1$.
\end{lemma}

\begin{proof}
 
Since $(\hat Q,\hat \delta)$ is a complex we must have 
$\pi_D(M_nf_{n-1} + f_n M_{n-2}) = 0$ for all $n>2$.  Applying $\phi$ we get
$\pi_A(M_n \phi(f_{n-1}) + \phi(f_n )M_{n-2}) = 0$.  

For $n>0$,  $\phi(\hat Q^n) = Q^n \oplus  Q^{n-1}(-1)$  and $\phi(\hat Q^0) = Q^0$.  The differential of $\phi(\hat Q)$, for $n\ge 2$, is  
$$ \phi(\hat \delta_n) = 
	\pi_A\left(\begin{array}{cc} M_n & \phi(f_n)\\ 0 & M_{n-1} \end{array}
	\right).
	$$
To ease notation, we supress $\pi_A$ in all of the following calculations. For $n\ge 2$, let $(p,p')$ be in the kernel of $\phi(\hat\delta_n)$.  Then 
$(pM_n,p\phi(f_n) + p'M_{n-1}) = (0,0)$.  Since $(Q,\delta)$ is acyclic, we may 
choose $q\in Q^{n+1}$ for which $qM_{n+1} = p$.  Then 
$(-q \phi(f_{n+1}) + p')M_{n-1} = qM_{n+1}\phi(f_n) + p'M_{n-1} = 0$ 
and we may choose $q'\in Q^n(-1)$ for which $q'M_n = (-q \phi(f_{n+1}) + p')$.  Clearly then 
$\hat \delta_{n+1}(q,q') = (p,p')$.   We have $H^n(\phi(Q)) = 0$ for all $n\ge 2$.  Since
$\phi(\hat \delta_1)$ is $\delta_1$ on $Q^1$ and zero on $Q^0(-1)$, we see that 
$H^0(\phi(\hat Q)) = K$.  

Finally,  
$ker(\phi(\hat\delta_1)) = ker(\delta_1)\oplus Q^0(-1) = im(\delta_2) \oplus A(-1)$.  
On the other hand,
$im(\phi(\hat \delta_2)) = im(\delta_2) \oplus (Q^2\phi(f_2) + Q^1(-1)M_1)$. 
By hypothesis, $Q^1M_1 = im(\delta_1) = A_+$.  By degree considerations, 
$Q^2\phi(f_2)$ must be contained in $Q^1(-1)M_1$.  Thus $H^1(\phi(\hat Q)) = K(-1)$.  
\end{proof}

It is clear from the proof of the Lemma that when $\hat Q$ is a complex, then $\phi(\hQ)$ is the mapping cone of a morphism of complexes $f:Q(-1) \to Q$. The following Lemma has the same proof as the Lemma above.

\begin{lemma}\label{complex2} Suppose that 
 $\hQ^3\to \hQ^2 \to \hQ^1 \to \hQ^0\to 0$ is a complex.  Then the complex
 $\phi(\hQ^3)\to \phi(\hQ^2) \to \phi(\hQ^1) \to \phi(\hQ^0)\to 0$ has homology:
$H^0 = K$, $H^1 = K(-1)$ and $H^2 = 0$.
\end{lemma}

It seems that the following lemma should be standard, but we include its proof.

\begin{lemma}\label{main1} If the central element $z$ is regular in $D$, then $(\hat Q,\hat \delta)$ is a complex and a minimal projective resolution of the trivial graded $D$-module $K$.
\end{lemma}

\begin{proof} Assume $z$ is regular.  We have
$$\begin{array}{lcl}z\cdot\pi_D(M_nf_{n-1} + f_n M_{n-2})&=&
 \pi_D(M_nzf_{n-1} + zf_n M_{n-2})\\	&=&
\pi_D( (-1)^{n}M_nM_{n-1}M_{n-2}  \\
&&\ \ \ \ \ \ \ + (-1)^{n-1}M_nM_{n-1}M_{n-2}) \\
 &=& 0
 \end{array}
 $$
Thus $ \pi_D(M_nf_{n-1} + f_n M_{n-2}) = 0$, and so by Remark \ref{corner}, 
$(\hat Q, \hat \delta)$ is a complex.  

By the definition of $\hat \delta_1$ we have $H^0(\hat Q) = K$.  Consider the short exact sequence of complexes 
 $$ 0 \to z \hat Q \to \hat Q \to \phi(\hat Q) \to 0.$$
Since $\hat Q$ is a complex of $z$-torsion-free modules, we have $z\hat Q \cong \hat Q(-1)$ and moreover $H(z \hat Q) = H(\hat Q(-1)) = H(\hat Q)(-1)$.  In particular, $H^0(z\hat Q) = K(-1)$.

If we apply this information to the long exact homology sequence induced by the short exact sequence above, recalling from Lemma \ref{complex1} that 
$H^0(\phi(\hat Q)) = K$,
$H^1(\phi(\hat Q)) = K(-1)$ and $H^n(\phi(\hat Q)) = 0$ for $n>1$, we obtain
$$0\to H^1(z\hQ) \to H^1(\hQ) \to K(-1) \to K(-1) \to K \to K \to 0$$
and for all $n>1$,
$$ 0\to H^n(z \hat Q) \to H^n(\hQ) \to 0. $$  Thus for all $n\ge 1$, 
$H^n(\hQ) = H^n(z \hQ) = H^n(\hQ)(-1)$.  Since all of these modules are bounded below in internal degree, this can only happen if they are all 0.  This proves that $(\hat Q,\hat \delta)$ is a projective resolution of $K$.  Since the modules in $\hat Q$ are all graded and every nonzero entry of the matrices $\hat M_n$ is homogeneous of degree at least 1, it follows immediately that the resolution is minimal.
\end{proof}

As we will see in the following theorem, the regularity of $z$ in $D$ can be detected entirely by the $\hat Q^3 \to \hat Q^2 \to \hat Q^1$ portion of $\hat Q$. 

\begin{thm}\label{main2}  Let $D$ be a central extension of $A$.  Then the following statements are equivalent:

(1) The central element $z\in D$ is regular in $D$.

(2)  The sequence of graded $D$-modules
$(\hat Q,\hat \delta)$ is a complex.  

(3)  The truncated sequence of graded $D$-modules
 $\hQ^3\to \hQ^2 \to \hQ^1 \to \hQ^0\to 0$ is a complex.  
 
(4) $\pi_D(M_3f_2 + f_3 M_1) = 0$ in $D$.

(5)  The sequence of graded $D$-modules $(\hat Q,\hat \delta)$ is a minimal projective resolution of $K$.

\end{thm}

\begin{proof}  In light of Lemma \ref{main1} and Remark \ref{corner}, we are required only to prove that (3) implies (1) and so we begin by assuming (3).  Let $\phQ$ be the truncated complex
 $\hQ^3\to \hQ^2 \to \hQ^1 \to \hQ^0\to 0$.

Let $L=Ann_D(z)$.  We must show that this graded two-sided ideal of $D$ is zero. Let 
$\alpha:D(-1) \to Dz$ be the epimorphism $b \mapsto bz$, with kernel $L(-1)$.  Since 
$\phQ$ consists of free $D$-modules, we can extend $\alpha$ appropriately and get a short exact sequence of complexes:
$$ 0 \to L\otimes_D \phQ(-1) \to \phQ(-1)\buildrel \alpha \over \to z\,\phQ \to 0.$$
We also have the short exact sequences of complexes
$$0 \to z\,\phQ \to \phQ \to \phi(\phQ) \to 0.$$ 
The long exact homology sequence of the second of these short exact sequences, taking into account Lemma \ref{complex2}, yields
$$ 0 \to H^1(z\, \phQ) \to H^1(\phQ) \to K(-1) \to H^0(z\, \phQ) \to K \to K \to 0.$$

Since $H^1(\phQ)= 0$ (by Remark \ref{corner}), we can now conclude 
$H^1(z\, \phQ) = 0$ and $H^0(z\, \phQ) = K(-1)$.

The long exact homology sequence of the first short exact sequence above now yields
$$0 \to H^0(L\otimes \phQ(-1)) \to H^0(\phQ(-1)) \to K(-1)\to 0.$$
Since $H^0('\hat Q(-1)) = H^0('\hat Q)(-1) = K(-1)$, we conclude
$H^0(L \otimes \phQ(-1)) = H^0(L\otimes \phQ)(-1) = 0$.  But 
$H^0(L\otimes \phQ) = L/D_+L$.  Since $L$ is bounded below, this can only be 0 if 
$L$ is 0, as required. 
\end{proof}

We have everything we need to prove Theorem \ref{main}, which follows from the restatement below.

\begin{thm}\label{mainagain} Let $A$ be a graded $K$-algebra of finite complexity 
$c=c(A)$.   Let $D$ be a central extension of $A$ such that the extending variable $z$ is not a zero divisor on $D_{\le c}$.  Then $z$ is regular in $D$.  
\end{thm} 

\begin{proof}   
From Theorem \ref{main2}, it suffices to prove 
 that  $\pi_D(M_3f_2 + f_3M_1)= 0 $.  But $z\cdot\pi_D(M_3f_2 + f_3M_1)= 0 $, as shown in the proof of \ref{main1}.  So it suffices to show that every entry 
 of  $M_3f_2 + f_3M_1$ has degree at most $c(A)$.  
 
Recall that $Q^n = A(m_{1,n},\ldots,m_{t_n,n})$.  We have
$m_{i,1} = -1$ for all $i$ and $\{-m_{1,2},-m_{2,2},\ldots,-m_{t_2,2}\}$ is the set of degrees of the relations $r_1,\ldots,r_m$ ($m=t_2$).  The $i$-$j$ entry of $M_2$ thus has degree $-1 - m_{i,2}$.  Similarly, the $i$-$j$ entry of $M_3$ has degree 
$m_{j,2}-m_{i,3}$.  It follows that the $i$-$j$ entry of $M_3M_2$ has degree $-1 - m_{i,3}$.   Since $\pi_D(M_3M_2 +f_3 z) = 0$ in $D$, the $i$-$j$ entry of $f_3$ must have degree $-2-m_{i,3}$.  Similarly, the $i$-$j$ entry of $f_2$ has degree $-1-m_{i,2}$.  Putting this all together, we see that the $i$-$j$ entry of $M_3f_2 + f_3M_1$ has degree $-1 -m_{i,3}$.  By hypothesis, this is at most $c(A)$, as required. 
\end{proof}

\section{Jacobi Conditions}

We return to the case of a deformation $U$ of the graded algebra $A$.  As before, let
$R=\{r_1,\ldots, r_m\}$ be the set of homogeneous relations of $A$ and $P=\{r_1+l_1,\ldots, r_m+l_m\}$ be the set of nonhomogeneous relations of $U$.  Let $D$ be the central extension associated to this deformation.  

We make an inductive definition of a nested family of subspaces of the free algebra $T$.  Let $V$ be the $K$-span of the generators $x_1,\ldots,x_n$ of $T$.  Let $P_1 = \hbox{span}_K(P \cap F^1T)$.  For $k>1$, set 
$$P_k = VP_{k-1} + P_{k-1}V + \hbox{span}_K(P \cap F^kT).$$

We see at once that $P_2$ is the span of the quadratic elements of $P$, $P_3$ the span of the quadratic and cubic elements of $P$ plus the cubic relations generated by $P_2$, etc.  It is clear that $P_1\subset P_2 \subset P_3\subset \ldots$ and that $\cup_k P_k = \langle P\rangle$.  It should also be clear that $P_k$ is not, in general, equal to $\langle P \rangle \cap F^kT$.  For example, if $P=\{x^2y-w, y^3 - w\}$, then $x^2w-wy^2 \in \langle P\rangle \cap F^3T$, but $P_3$ is just the $K$-span of $P$, which does not contain $x^2w-wy^2$. 

We need one simple Lemma.  Let $\phi_1:T[z] \to T$ be the epimorphism given by evaluation at $z=1$.  Recall the homogenization function $h:T\to T[z]$.  We omit the proof.

\begin{lemma}\label{jaclemma}  For all $k\ge 1$, $\phi_1(\langle h(P)\rangle_k) = P_k$.  
\end{lemma}

We can now restate Theorem \ref{main} in the language of a Jacobi type condition.

\begin{thm}\label{jacobi}  Let $A$ be a graded $K$-algebra of finite complexity $c(A)$ and let 
$U$ be a deformation of $A$.    Then $U$ is a PBW-deformation of $A$ if and only if 
$P_1=0$ and the following Jacobi condition is satisfied:
$$ P_{k+1} \cap F^kT \subset P_{k} \ \ \hbox{ for all }\ \ 1\le k \le c(A).$$
\end{thm}

\begin{proof}
Let $D$ be the central extension of $A$ associated to the deformation $U$.  By combining theorems \ref{main} and \ref{standard}, we see that we are required to prove that $z$ is not a zero divisor on $D_{\le c(A)}$ if and only if the Jacobi condition is satisfied.  

Assume the Jacobi condition and suppose also that $z$ is a zero divisor on $D_{\le c(A)}$.  Then for some $x \in T[z]_k$, $k\le c(A)$, with $x\not\in\langle h(P)\rangle$, we have $zx \in \langle h(P)\rangle_{k+1}$.  By lemma \ref{jaclemma}, we have $\phi_1(zx) = \phi_1(x) \in P_{k+1} \cap F^kT$.  Therefore
$\phi_1(x)\in P_k$.  Again by \ref{jaclemma}, we can choose 
$y\in \langle hP\rangle_{k}$ with $\phi_1(y)=\phi_1(x)$.  But the restriction of $\phi_1$ to any homogeneous component of $T[z]$ is injective.  Thus $x=y$, a contradiction.

Conversely, assume the Jacobi condition fails. Then for some $k\le c(A)$,  there is a 
$w \in P_{k+1}\cap F^kT$ such that $w\not\in P_k$.  Let $x=h(w)\in T[z]_k$.  Then $\phi_1(x) = w$, so $x\not\in \langle h(P)\rangle_k$, but $\phi(zx) = w \in P_{k+1}$, so
$zx \in \langle h(P)\rangle_{k+1}$.  Thus $z$ is a zero divisor on $D_{\le c(A)}$.
\end{proof}

\begin{rmk}\label{Jacobi1} 
If $A$ is $N$-Koszul then  $P_i=0$  for $i<N$ and $P_N$ is the $K$-span of $P$.   It then follows directly from Theorem \ref{jacobi}  that the Jacobi condition given in \cite{BGZ} and \cite{FV} is equivalent to $z$ being regular in degree in $N$.   Hence the PBW-deformation theorems in \cite{BGZ} and \cite{FV}  for $N$-Koszul algebras can be recovered using the regularity theorem \ref{main}.
\end{rmk}

\section{Examples}

We conclude with a very short list of examples. 

 \begin{ex}
Three dimensional Artin-Schelter (AS) regular algebras \cite{AS} are either quadratic or cubic.  The quadratic ones are Koszul and hence have complexity two, while the cubic ones are 3-Koszul and  have complexity three.  

Four dimensional AS regular algebras, on the other hand, may have two, three or four generators.  The 4-generated algebras are Koszul.  The 3-generated algebras have two relations of degree 2 and 2 of degree 3 and satisfy a resolution of the form
$$0\to A(-5) \to A(-4)^3 \to A(-3,-3,-2,-2)\to A(-1)^3\to A\to K\to 0.$$
The complexity of such algebras is thus 3. Perhaps the most interesting case, from the complexity point of view, is the 2-generated AS regular algebras of dimension 4, which have one cubic relation, and one quartic relation.  These algebras have been studied extensively by Lu et. al. \cite{LPWZ}.  For such  an algebra there is a minimal resolution of the form
$$ 0\to A(-7)   \to A(-6)^2\to A(-3,-4)\to A(-1)^2  \to A \onto K. $$
Hence these algebras have complexity 5.

\end{ex}
 
Infinite complexity is easily arranged.

 \begin{ex}  Let $A=K\la w,x,y,z\ra / \la yz, zx-xz,zw\ra$.  Using the facts that $x^nz\ne 0$ in $A$ but  
 $yx^nz=0$ in $A$, one can establish that $A$ has global dimension 3 and a minimal projective resolution:
 $$ 0 \to A(-3,-4,-5,\ldots ) \to A(-2)^3 \to A(-1)^4 \to A \to K \to 0.$$
 Hence $c(A)$ is infinite.  
  \end{ex}
 
 And finally, using an essentially trivial example, we can show that even when the complexity of $A$ is not 0, it need not be as large as the maximal degree of the relations of $A$.  
 
 \begin{ex}  Let $A = K\la x,y,z \ra/\la y^2,xyz \ra$.  Then $A$ has relations of degrees 2 and 3, but the minimal resolution associated to $A$ is
 $$ \ldots A(-5) \to A(-4) \to A(-3) \to A(-2,-3) \to A(-1)^3 \to A \to K \to 0$$
 and so the complexity of $A$ is only 2. 
 \end{ex}


\bibliographystyle{amsplain}

\bibliography{bibliog}

\end{document}